\documentclass[12pt]{article}
\usepackage{epsfig}
\usepackage{amsfonts}
\newtheorem{theorem}{Theorem}[section]
\newtheorem{lemma}[theorem]{Lemma}

\newtheorem{proposition}[theorem]{Proposition}
\newtheorem{corollary}[theorem]{Corollary}

\newcommand{\R}{{\mathbb R}}

\newcommand{\Z}{{\mathbb Z}}
\newcommand{\N}{{\mathbb N}}
\newcommand{\be}{\begin{equation}}
\newcommand{\ee}{\end{equation}}

\newcommand\vnz{V_{N,\underline{Z}}}
\newcommand\an{{\cal A}_N}

\newcommand\munz{\mu^{(N,\underline{Z})}}
\newcommand\munjzj{\mu^{(N_j,\underline{Z}_j)}}
\newcommand\munznu{\mu^{(N_\nu,\underline{Z}_\nu)}}
\newcommand\xonenz{x_1^{(N,\underline{Z})}}
\newcommand\xnnz{x_N^{(N,\underline{Z})}}
\newcommand\weakstar{\rightharpoonup^*}
\newcommand\M{{\cal M}}
\newcommand\inz{\tilde{V}_{(N,\underline{Z})}}
\newcommand\Ilambdaz{I_{\lambda,\underline{z}}}
\newcommand\inznu{I^{(N^{(\nu)},\underline{Z}^{(\nu)})}}

\newcommand\phieps{{\phi_\varepsilon}}
\newcommand\mueps{{\mu^\varepsilon}}
\newcommand\Mlambda{{\cal M}_{[0,\lambda]}}
\newcommand\Qi{Q^{(i)}}

\newcommand\diam{{\mbox{diam}\,}}
\newcommand\supp{{\mbox{supp}\,}}
\newcommand\Itilde{\tilde{I}}

\newcommand\xik{x_i^{(k)}}

\newcommand\Diag{\mbox{diag}}

\title{Minimum energy configurations of classical charges: \\ Large $N$
asymptotics}
\author{ {\large Stephane Capet and Gero Friesecke}  \\[2mm] 
Center for Mathematics, Technische Universit\"at M\"unchen, \\ 85747 Garching, Germany}
\date{June 10, 2009}
\begin{document}
\maketitle

\begin{abstract}
We study minimum energy configurations of $N$ particles in $\R^3$ of charge $-1$ (`electrons')
in the potential of $M$ particles of charges $Z_\alpha>0$ (`atomic nuclei'). 
In a suitable large-N limit, we determine the asymptotic electron
distribution explicitly, showing in particular that the number of 
electrons surrounding each nucleus is asymptotic to the nuclear charge (``screening'').
The proof proceeds by establishing, via Gamma-convergence, a coarse-grained variational 
principle for the limit distribution, which can be solved explicitly. 
\end{abstract}

\section{Introduction}
The goal of this paper is to shed new light on basic screening effects
in molecules. By screening one means the remarkable tendency of electrons to usually group themselves
around the atomic nuclei in such a way so as to cancel much of the long range $\sim 1/R$ Coulomb potential exerted
by the nuclei and make the net potential exerted by the atoms short-range. 

Screening is usually tacitly assumed in molecular mechanics, molecular dynamics, 
statistical mechanics, and continuum mechanics. One starts from the outset from
short-range atomistic forces respectively short-range
continuum forces (i.e. stresses alias surface forces). 

Large-scale failure of screening
(which is not observed in nature) 
would lead to spectacular breakdown of these models, e.g. bare Coulomb interactions 
violate the linear scaling of energy $E$ with volume $V$ for a quantum mechanical
crystal. To see this, place bare atomic nuclei on the integer lattice points in a 3D cube
of sidelength $L$, $\{R_1,...,R_M\}=\Z^3\cap [0,L]^3$ and evaluate their interaction
energy asymptotically in the limit of large $L$: 
\begin{eqnarray}
   E & \!\sim\! & \!\mbox{$\sum_{{R_\alpha,\, R_\beta\in \Z^3}\atop{|R_\alpha|,\, |R_\beta|\le L}}  
                        \!\frac{1}{|R_\alpha-R_\beta|}$}  \sim \mbox{$
                        \int\!\!\int_{|x|,|y|\le L}\frac{1}{|x-y|}d(x\!-\!y)d(x\!+\!y)$} \nonumber \\
               & \!\sim\! & L^3 \cdot L^2 \; \sim \; M^{5/3} \; \sim \;  V^{5/3}, \label{catastrophe}
\end{eqnarray}
i.e. the energy per atom tends to infinity as the system gets large (in a finite system
of $M=10^{23}$ atoms it is already too large by a factor of about $10^{17}$).

The above example is unstable, but small-scale failure of screening is common in nature, and
yields important $O(1)$ contributions to the energy per atom.
Examples include ionic crystals like
NaCl, molecules with low permanent multipole moment like H$_2$O, intermediate states during 
chemical reactions, and core regions of atoms. 

We know of no mathematical results which directly explain and quantify screening from full quantum mechanics. 
The perhaps furthest result in this direction concerns an indirect, coarse-scale manifestation of screening:
the ground state energy of a molecule with $M$ atoms is known not to scale like the example (\ref{catastrophe}),
but is bounded above and below by a constant times $M$ [DL67, LD68, LT75]. 

More insight has been obtained in asymptotic limits. For atoms in the limit
of large atomic number $Z$, it is known [LS77] that the total electron density
is asymptotically radial, the profile given by Thomas-Fermi theory, and falls off like
$r^{-6}$, and so the asymptotic net potential exerted by the atom is short-range. For
a closely related result see [ILS96]. Another interesting limit is the thermodynamic limit
for crystalline solids, in which the nuclei are arranged
on a regular subset of a periodic crystal lattice, say $B_R\cap\Z^3$, where $B_R$ denotes the ball of radius $R$ around
the origin, and $R$ tends to infinity. In this case, for a slightly simplified version of quantum mechanics (absence of spin, a rigid wall
assumption, and coupling to an electron reservoir), the ground state energy is known to be asymptotically proportional
to the number of nuclei [Fe85]. Moreover for convex density functional models such as the Thomas-Fermi-Weizs\"acker model, the
ground state density is known to become asymptotically periodic [CLL98]. The latter result, by Catto, Le Bris and
Lions, can be viewed as a quantitative version of screening: the asymptotic amount of electron density in each unit cell exactly cancels the nuclear charge in that cell,
making the net electrostatic potential excerted by the cell short-range. 

Here we introduce and analyze a model which allows some new insight into screening for general, non-periodic, arrangements
of nuclei, at the expense of further simplification of the treatment of electrons.
The model maintains the long-range Coulomb forces between electrons and
atomic nuclei exactly, but treats the electrons as classical point charges and
replaces the Laplacian in the electronic Schr\"odinger equation by a hard-core constraint.
Our main results are 
\begin{itemize}
\item 
a simple proof of exact screening in a large-$N$ continuum limit of this model, via explicit
determination of the minimizer (see (\ref{explicit}))
\item
a proof of approximate screening in this model for large finite $N$. This is done by establishing, via Gamma-convergence, 
that the discrete minimizers converge to the minimizer of the continuum limit (see Theorems \ref{thm1},
\ref{thm2}).
\end{itemize}
Our model is variational, and describes a system of $N$
particles in $\R^3$ of charge $-1$ (`electrons'), with variable positions $x_1,...,x_N$, 
which Coulomb-repel each other and are 
Coulomb-attracted to M particles of charges $+Z_\alpha$ at fixed positions $R_\alpha$ 
(`atomic nuclei'): 

Minimize
\be \label{vnz}
   \vnz(x_1,..,x_N) := \sum_{i=1}^N v(x_i) + \sum_{1\le i<j\le N} 
   \frac{1}{|x_i-x_j|}, \hspace*{1cm}(N\in\N)
\ee
where
\be \label{molpot}
    v(x) = \sum_{\alpha=1}^M \frac{-Z_\alpha}{|x-R_\alpha|} \;\;\;\;\;\;\;\;\;\; 
    (\underline{Z}=(Z_1,..,Z_M), \; Z_\alpha>0, \; R_\alpha\in\R^3),
\ee
over the set
\be \label{an}
   \an := \{ (x_1,..,x_N)\in\R^{3N} \, \Bigl|\Bigr. \, |x_i-R_\alpha|\ge d \mbox{ for all
   }i,\,\alpha \} \;\;\;\;\; (d>0).
\ee
The hard core
assumption (\ref{an}) may be viewed as a 
crude ``uncertainty principle'' which prevents
electrons from falling into the nucleus, with the hard core radius $d$ playing the
role of $\hbar$. More precisely, the model
(\ref{vnz}), (\ref{molpot}), (\ref{an}) arises from the full 
quantum mechanical (Born-Oppenheimer-)Hamiltonian of the electrons in a molecule,
\be \label{ham}
   H_{N,\underline{Z}} = -\frac{1}{2}\Delta + \vnz 
           = \sum_{i=1}^N\Bigl(-\frac12\Delta_{x_i} +v(x_i)\Bigr) + \sum_{1\le i<j\le N} 
   \frac{1}{|x_i-x_j|},
\ee
by replacing the one-body operator
$-\frac12\Delta_{x_i}+v(x_i)$ by the effective potential $v_{eff}(x_i):= v(x_i)$ when $|x_i-R_\alpha|\ge d$
for all $i$ and all $\alpha$,
$+\infty$ otherwise. 

The physics of the model (1), (2), (3) is independent of the choice of hard core radius $d$, as long as the hard cores
are not overlapping, i.e. 
\be \label{nonoverl} d\in\left\{\begin{array}{ll}(0,\infty) & \mbox{if }M=1\\
      (0,\mbox{$\frac12$}\min_{\alpha\neq\beta}|R_\alpha-R_\beta|) & \mbox{if }M\ge 2.\end{array}\right.
\ee
A different choice
just corresponds to an overall scale factor of length and energy.

We proceed to describe our results on the model (2), (3), (4). Comparisons with what is known (or expected) 
in quantum mechanics and related models are postponed to the end of this Introduction.  

Our first two observations, the second of which is at first sight somewhat surprising, are the following:
\begin{proposition} \label{prop1} a) (Attainment for neutral molecules and singly-negative ions) Let $Z:=\sum_{\alpha=1}^M
Z_\alpha$.
For $N\le Z+1$, there exists a minimizer of $\vnz$ on $\an$. \\[2mm]
b) (Absorption principle) 
Every minimizer $(x_1,..,x_N)$ of $\vnz$ on $\an$ satisfies $x_i\in\bigcup_{\alpha=1}^MS_\alpha$ for all $i$,
where $S_\alpha=\{x\in\R^3 \, | \, |x-R_\alpha|=d\}$ denotes the sphere of radius $d$ centred at $R_\alpha$.
\end{proposition}
\setcounter{theorem}{0}
{\bf Proof} a) follows from standard arguments in the calculus of variations. The fact that unlike in quantum
mechanics, attainment can also be shown for $N=Z+1$ comes from the fact that the joint potential exerted by
$Z$ particles at $x_1,..,x_Z$ and the nucleus onto an additional particle on a sphere of radius $R>\max |x_i|$ 
is zero on average but nonconstant, and hence negative somewhere. We omit the
details. b) follows from observing that the potential
$\vnz(x_1,...,x_N)$ is a harmonic function with respect to each particle position $x_i$ and applying the maximum
principle. See e.g. [Lan72] for a related observation for purely repulsive Coulomb particles confined to a bounded (instead
of unbounded) set.
\\[2mm]
b) allows to allocate each electron unambiguously to one atomic nucleus. Investigating the extent of screening
means investigating \\
--- how closely the number of electrons going to any particular nucleus matches the nuclear
charge (a perfect match corresponds to a zero net monopole moment of the atom) \\
--- how uniformly and symmetrically the electrons distribute themselves around the nucleus (this determines
the higher net multipole moments).

In the asymptotic limit when the nuclear charges are large, these questions
have simple answers.
\begin{theorem} \label{thm1} Assume $d$ satisfies (\ref{nonoverl}), and denote $Z:=\sum_{\alpha=1}^MZ_\alpha$.
Let $(x_1^{(N,\underline{Z})},..,x_N^{(N.\underline{Z})})$ be any minimizer of $\vnz$. In the limit
\be \label{nlimit}
   N=Z\to\infty, \;\;\; \frac{Z_\alpha}{Z} \to z_\alpha,
\ee
a) (Neutrality) 
$$
    \frac{\sharp \{ x_i^{(N,\underline{Z})}\, | \, x_i^{(N,\underline{Z})}\in S_\alpha\} }
    {Z_\alpha} \longrightarrow 1
$$
b) (Equidistribution) For any measurable $\Omega\subseteq S_\alpha$
with $area(\partial\Omega)=0$, 
$$
    \frac{ \sharp \{ x_i^{(N,\underline{Z})}\, | \, x_i^{(N,\underline{Z})}\in\Omega\} } {N} \longrightarrow z_\alpha
    \frac{\mbox{area}\, (\Omega)}{\mbox{area}\, (S_\alpha)}
$$
c) (Limit energy)  
$$
    \frac{\vnz(\xonenz,..,\xnnz)}{N^2} \longrightarrow -\frac{1}{2}\sum_{\alpha=1}^M\frac{z_\alpha^2}{d} 
    - \sum_{1\le \alpha<\beta\le M}\frac{z_\alpha z_\beta}{|R_\alpha-R_\beta|}.
$$
\end{theorem}
Physically, the results in a), b) and c) 
are ``screening results'' which 
mean, respectively, that in the above limit\\
-- the monopole moment of each atom vanishes \\
-- the higher multipole moments of each atom vanish \\
-- the interaction energy between the atoms vanishes. \\
To understand this interpretation of c), consider,
instead of the electronic energy $\vnz$, the total classical energy of the molecule which includes the Coulomb
repulsion between the nuclei, 
\be \label{Eclass} 
    E^{class}_{N,\underline{Z}} = \inf_{\an}\vnz + V^{nuc}_{\underline{Z}}, \;\;\;\;\;
        V^{nuc}_{\underline{Z}}=\sum_{1\le\alpha<\beta\le M}\frac{Z_\alpha
        Z_\beta}{|R_\alpha-R_\beta|}.
\ee
(Here and  below we use the convention that $V^{nuc}_{\underline{Z}}=0$ when $M=1$.)
The formula in c) then says that 
$$
    \frac{E^{class}_{N,\underline{Z}}}{N^2} \longrightarrow -\frac12\sum_{\alpha=1}^M \frac{z_\alpha^2}{d}.
$$
(This is because by (\ref{nlimit}), $Z_\alpha Z_\beta/N^2$
converges to $z_\alpha z_\beta$ and hence the internuclear repulsion term cancels the second term appearing in c).)
%
%
%
%
%
%
In other words, the limit energy of the molecule just equals the sum of the limit 
energies of the individual atoms. In particular, it is independent of the atomic positions $R_1$,...,$R_M$,
i.e. contains no interaction terms. 

The above screening results are a corollary of the following more general result, which in
addition uncovers interesting behaviour of excess charges moving off to infinity in
case of negative ions $N>Z$. 

To include this case we consider, instead of $N=Z\to\infty$, the more general limit
\be \label{limit}
   N\to\infty, \;\;\; Z=\sum_{\alpha=1}^MZ_\alpha\to\infty, \;\;\; \frac{N}{Z} \to\lambda, \;\;\; \frac{Z_\alpha}{Z} \to z_\alpha,
\ee
where $\lambda\in(0,\infty)$ is a filling factor. Positive ions correspond to $\lambda<1$, neutral molecules
to $\lambda=1$, and negative ions to $\lambda>1$. 

For negative ions, the minimum of $\vnz$ on $\an$ is typically not attained (see below) and so
one needs to relax the restriction to exact minimizers in Theorem \ref{thm1}. Instead one considers
more general low-energy states, in the sense of 
$$
   \mbox{energy difference from infimum} \; < < \; \mbox{total energy},
$$
as made precise by the following \\[2mm]
{\bf Definition:} A sequence $\{(\xonenz,..,\xnnz)\}$ is called
a sequence of approximate minimizers of $\vnz$ in the limit (\ref{limit}) if
\be \label{approxmin}
    \frac{ \vnz(\xonenz,..,\xnnz)-\inf_{\an}\vnz } {Z^2} \longrightarrow 0.
\ee
\begin{theorem} \label{thm2} (Variatonal principle for the limit distribution)
For any sequence $\{(\xonenz,..,\xnnz)\}$ of approximate minimizers of $\vnz$, 
in the limit (\ref{limit}) the associated measures
\be \label{measures}
   \munz := \frac{1}{Z} \sum_{i=1}^N \delta_{x_i^{(N,\underline{Z})}}
\ee
satisfy 
\be\label{measconv}
   \munz \weakstar \mu_\lambda
\ee
and
\be
   \frac{\vnz(\xonenz,..,\xnnz) } {Z^2} \longrightarrow  \Ilambdaz(\mu_\lambda),
\ee 
where $\Ilambdaz\, : \, {\cal M}_+(\R^3\backslash \Omega)\to \R\cup\{+\infty\}$ (see below for notation)
is the continuum energy functional 
\be \label{contlimit}
   \Ilambdaz(\mu) := \left\{ \begin{array}{ll} -\int_{\R^3\backslash \Omega} \sum_{\alpha=1}^M \frac{z_\alpha}{|x-R_\alpha|} d\mu(x) 
               +\frac{1}{2}\int\int_{(\R^3\backslash \Omega)^2 }
               \frac{1}{|x-y|}d\mu(x)\, d\mu(y) & \mbox{if }\int d\mu\le\lambda, \\
              +\infty & \mbox{otherwise}, \end{array}\right. 
\ee
and $\mu_\lambda$ is its unique minimizer. 
\end{theorem}
Existence of a unique minimizer of $\Ilambdaz$ is proved in Proposition 2.1 below. Here and
below our notation is as follows: the halfarrow $\weakstar$ denotes weak* convergence in the 
space $\M(\R^3\backslash \Omega)$
of Radon measures on $\R^3\backslash \Omega$,\footnote{Recall that
for any closed subset $A\subseteq\R^d$, $\M(A)$ 
is the dual of the space $C_0(A)=\{f\, : \, A\to\R\, | \, 
f \mbox{ continuous},
\,f(x)\to 0\mbox{ for }|x|\to\infty\}$, and that a sequence of Radon measures $\mu_\nu$ is said to
converge weak* to $\mu$, notation: $\mu_\nu\weakstar\mu$, if 
$\int_{A}f\, d\mu_\nu \to \int_{A}f\, d\mu$ for all $f\in C_0(A)$.}
${\cal M}_+(\R^3\backslash\Omega)$ denotes the set $\{\mu\in {\cal M}(\R^3\backslash\Omega) \, | \, 
\mu\ge 0$, $\Omega$ is the union of the hard cores of the
nuclei, i.e. $\Omega=\cup_{\alpha=1}^MB_d(R_\alpha)$, 
$B_d(R_\alpha)=\{x\in\R^3\, | \, |x-R_\alpha|<d\}$
and $\underline{z}=(z_1,...,z_M)$.

The point about Theorem 1.2 is that the electrostatic continuum energy $\Ilambdaz$ which appears in the limit
is much simpler than the intricate particle energy $\vnz$. For neutral
molecules $(\lambda=1)$ or negative ions $(\lambda>1)$, and nonoverlapping hard cores (i.e. (\ref{nonoverl})),
the minimizer of the continuum energy can be determined explicitly, 
\be \label{explicit}
   \mu_\lambda \equiv \mu_1 = \sum_{\alpha=1}^M z_\alpha \frac{H^2\bigl|_{S_\alpha}\bigr.}{4\pi d^2},
\ee
where $H^2|_{S_\alpha}$ denotes two-dimensional Hausdorff measure restricted to the sphere $S_\alpha=\{x\in\R^3 \, | \, 
|x-R_\alpha|=d\}$
(see Proposition 2.1 e) below). This together with (\ref{measconv}) readily implies the screening results
in Theorem 1.1 (see Section 2).

Theorem 1.2 together with formula (\ref{explicit}) also lead to interesting conclusions about instability
of negative ions. Note that for negative ions ($\lambda>1$) the limit measure has less mass
than the approximating measures,
$$
   \int d\mu_\lambda = 1 < \lambda = \lim\frac{N}{Z} = \lim \int d\munz.
$$
(The first equality is due to the fact that $\sum_{\alpha=1}^Mz_\alpha = \lim \sum_{\alpha=1}^M\frac{Z_\alpha}{Z} = 1.$)
Physically this means that only $Z+o(Z)$ particles stay bound and $N-(Z+o(Z))$ particles move off to infinity.
For a precise formulation (as a nonattainment theorem for $\vnz$ when $Z$ is sufficiently large and $N$ exceeds $Z$
by a nonzero fraction) see Section 5.

We establish Theorem \ref{thm2} by showing that the particle energy $\vnz$ and the continuum energy $\Ilambdaz$
are related in the mathematically rigorous sense
of Gamma-convergence, introduced by De Giorgi (see [DM88, Br02] or the beginning of Section 5).
Starting point is the observation that the particle energy $\vnz$ can be re-interpreted in a natural
way as an energy functional on the space ${\cal M}_+(\R^3\backslash \Omega)$ of nonnegative Radon measures on
$\R^3\backslash \Omega$ (with $\Omega$ as defined below Theorem 1.2). Define
\be \label{continuummodel}
   \inz(\mu) := 
    - \int_{\R^3\backslash \Omega} \sum_{\alpha=1}^M \frac{Z_\alpha}{Z}\frac{1}{|x-R_\alpha|} d\mu(x) 
               +\frac{1}{2}\int\int_{(\R^3\backslash \Omega)^2\backslash\mbox{diag} }
               \frac{1}{|x-y|}d\mu(x)\, d\mu(y)
\ee
if $\mu=\frac{1}{Z}\sum_{i=1}^N\delta_{x_i}$ for some distinct $x_1,..,x_N\in\R^3\backslash\Omega$,
and set $\inz(\mu) := +\infty$ otherwise. Here and below diag denotes the diagonal 
$\{(x,x)\, | \, x\in\R^3\backslash \Omega\}$. Then for $\mu$ as in the first 
alternative, we have the identity 
\be \label{Inzid}
   \inz(\mu) = \frac{1}{Z^2}\vnz(x_1,...,x_N).
\ee
We then show:
\begin{theorem} \label{theo3} (Gamma-convergence)
In the limit (\ref{limit}), the sequence of functionals $\inz\, : \, {\cal M}_+(\R^3\backslash \Omega)\to \R\cup\{+\infty\}$ 
Gamma-converges (with respect to weak* convergence of Radon measures) 
to the functional $\Ilambdaz$ defined in (\ref{contlimit}). 
\end{theorem}
Physically, Theorem \ref{theo3} means that the limit functional not just correctly captures the ground state energy,
but also any energy change by a nonvanishing fraction of the ground state energy when the ground state is deformed.

Note that the restriction to discrete measures has disappeared, and
in the domain of integration of the second term the diagonal is now included.
The latter is essential, for otherwise the functional would promote clustering rather than equidistribution,
and the minimizers would, e.g. in the case of atoms $(M=1)$, be given by $\lambda\delta_x$, where $x$ is any point
on the sphere $|x-R_1|=d$. In particular an unlimited amount of electronic charge could
be bound by the nucleus.

We proceed to compare our results to various results in the literature on other models with many-body Coulomb interactions.
The attainment result of Proposition \ref{prop1}a) continues to hold in quantum mechanics
(Zhislin's theorem, see e.g. [Fr03] or see the original Russian article [Zh60]), but requires the slightly stronger hypothesis $N<Z+1$ which excludes
singly-negative ions. The ``no shells'' result of Proposition \ref{prop1}b) is false for true atoms (see [BB55] for
experimental data showing multiple maxima of the radial electron density in Argon), but
interestingly, it is also false, e.g., for classical Coulomb particles confined to a disc in two dimensions, in which case
minimizers would extend into the radial direction [EO00]; but it would become true again if the interaction was replaced by the
Green's function of the {\it two}-dimensional Laplacian. 
For a result related to Theorem \ref{thm2} for 
repulsive classical charges confined to a compact set see [Lan72], where it is proved that every sequence of empirical
measures of minimizers of the particle system contains a subsequence converging to a minimizer of the relevant continuum limit.
We know of no analogues, neither classical nor quantum, of the Gamma-convergence result of Theorem \ref{theo3}.
The nonattainment result of Corollary \ref{cor1} for $Z$ large and $N>Z+o(Z)$ (see the discussion following Theorem \ref{thm2}) 
is known to hold analogously in quantum mechanics in the special case $M=1$ [LSST], the case $M>1$ being open.
Numerical data of [MDH96] for $N$-particle configurations with minimal Coulomb repulsion on the sphere suggest that the precise
attainment threshold $N(Z)$ of our classical model with $M=1$ equals $Z$ plus a slowly growing function of $Z$. For an
investigation of the higher order energy asymptotics of the latter problem see [KS98]. 

Finally we remark that the passage from the particle energy (\ref{vnz}) to the continuum energy (\ref{continuummodel}) 
is achieved here via a conceptually new viewpoint which should be of more general interest. Instead of parametrizing particle configurations
$(x_1,..,x_N)$ with respect to some reference configuration (Lagrangian viewpoint), one considers the associated empirical measure 
$const\sum_{i=1}^N\delta_{x_i}$ which counts how many particles are contained in a given spatial region (Eulerian viewpoint). 
This allows us here to pass to the continuum limit for
a frame-indifferent particle system with re-labelling symmetry, without any a priori assumptions on
admissible particle configurations. This strategy should be applicable, at least in principle,
to other interesting problems, such as many-atom systems interacting via Lennard-Jones-type potentials. For why the latter problem,
despite involving short-range rather than long-range interactions, is in fact harder, see the remarks at the beginning of Section 5.2. 

Our plan of the paper is as follows. In Section 2 we analyze the limit theory.
In Section 3 we show how the asymptotic results in theorems \ref{thm1}
and \ref{thm2} follow from Gamma-convergence and the explicit solution of the limit theory. Section 4 is devoted to proving 
nonattainment for sufficiently negative ions, and in the final section we establish the Gamma-convergence result of Theorem \ref{theo3}.

\section{Analysis of the limit theory}
Here we analyze the limit theory (postponing its rigorous justification to Sections 4 and 6). 
Its fundamental advantage over the particle system is that it can be minimized
explicitly. 

\begin{proposition} \label{limitmodel} (Analysis of limit theory) For any $\lambda\ge 0$, $z_1,..,z_M\ge 0$, $R_1,..,R_M\in\R^3$, $d>0$,
and with $\Omega$ as defined below Theorem \ref{thm2}, the functional 
$$
   I(\mu) = -\int_{\R^3\backslash \Omega} \sum_{\alpha=1}^M\frac{z_\alpha}{|x-R_\alpha|}d\mu(x) 
            + \frac12 \int\int_{(\R^3\backslash \Omega)^2} \frac{1}{|x-y|}d\mu(x)d\mu(y)
$$
possesses a unique minimizer $\mu_\lambda$ on $\Mlambda = \{\mu\in {\cal M}(\R^3\backslash \Omega)\, | \, \mu\ge 0, \int d\mu \le \lambda\}$.
Moreover when condition (\ref{nonoverl}) holds, the following statements are true: 
\begin{itemize}
\item[a)] (saturation of mass) $\int d\mu_\lambda = \min\{\lambda,z\}$, where $z:=\sum_{\alpha=1}^Mz_\alpha$.
\item[b)] (saturation of minimizer) $\mu_\lambda=\mu_z$ for all $\lambda\ge z$. 
\item[c)] (saturation of energy) $e(\lambda):=\min_{\Mlambda}I$ is continuous, strictly decreasing for $\lambda\in[0,z]$, and constant for
$\lambda\in[z,\infty)$. 
\item[d)] (Minimizer for atoms) If $M=1$, 
$$
   \mu_\lambda = c(\lambda) \frac{H^2\bigl|_{S}\bigr.}{4\pi d^2}, \;\;\; e(\lambda) = -\frac{zc(\lambda)}{d} + \frac{c(\lambda)^2}{2d},
$$
where $c(\lambda)=\min\{\lambda,\, z\}$ and $S$ is the sphere of radius $d$ centred at $R_1$.
\item[e)] (Minimizer for molecules) If $M\ge 2$, $\lambda\ge z$ (neutral molecules or negative ions), 
$$
   \mu_\lambda \equiv \mu_z = \sum_{\alpha=1}^M z_\alpha \frac{H^2\bigl|_{S_\alpha}\bigr.}{4\pi d^2}, 
   \;\;\; e(\lambda)\equiv e(z) = - \sum_{\alpha=1}^M\frac{z_\alpha^2}{2d} - \sum_{1\le\alpha<\beta\le M} 
    \frac{z_\alpha z_\beta}{|R_\alpha-R_\beta|},
$$
where $S_\alpha$ is the sphere of radius $d$ centred at $R_\alpha$. 
\end{itemize}
\end{proposition}
Note that $I$ is well-defined on ${\cal M}_+(\R^3\backslash\Omega)$ 
(the space of nonnegative Radon measures of finite mass on $\R^3\backslash \Omega$, as introduced below Theorem \ref{thm2}) as an element
of $\R\cup\{+\infty\}$, because the negative term $-\sum_\alpha\int_{\R^3\backslash \Omega} z_\alpha|x-R_\alpha|^{-1}d\mu(x)$ is 
always finite, due to the boundedness of the integrand on the domain of integration.

The proof of the proposition relies on developing some basic functional analysis for the Coulomb self-energy
functional 
\be \label{Jfunctional}
   J(\mu) := \frac{1}{2}\int\int_{\R^6}\frac{1}{|x-y|}d\mu(x)d\mu(y)
\ee
on Radon measures. For smooth, rapidly decaying measures, this functional reduces to the very 
well studied Dirichlet integral for the associated potentials. But this standard setting is insufficient here,
as the minimizers themselves are singular measures which concentrate on lower-dimensional surfaces (see d) and e)).
This reflects the fact that we are dealing with a support constraint on the measures, as opposed to,
say, boundary conditions on the potential.

Denote by ${\cal C}(R^3)$ the set of nonnegative
Radon measures on $\R^3$ of finite mass for which $J(\mu)$ is finite.
Define an extension of $J$ to measures with both negative and positive
part, as follows: if $\mu=\mu_1-\mu_2$ with $\mu_1$, $\mu_2\in{\cal C}(R^3)$, set
\begin{eqnarray*}
  J(\mu_1-\mu_2) & := & \frac12 \int\int_{\R^6}\frac{1}{|x-y|}d\mu_1(x)d\mu_1(y)
                    -  \int\int_{\R^6}\frac{1}{|x-y|}d\mu_1(x)d\mu_2(y) \\
                 &   + & \int\frac12 \int_{\R^6}\frac{1}{|x-y|}d\mu_2(x)d\mu_2(y).
\end{eqnarray*}
Since the first and last term are finite by assumption, and the integrand in the middle term is nonnegative,
this is well defined as an element of $\R\cup\{-\infty\}$. The key property of $J$ needed in the proof of the proposition is
\setcounter{theorem}{0}
\begin{lemma} \label{positive} $J(\mu_1-\mu_2)\ge 0$ for any $\mu_1$, $\mu_2\in{\cal C}(\R^3)$, with equality
if and only if $\mu_1=\mu_2$. 
\end{lemma}
This result is trivial for smooth, rapidly decaying measures, as well as for nonnegative measures. That it should continue
to hold for rough measures without a sign is well known ``folklore'' in part of the potential theoretic literature
(see e.g. [La75]). Our proof, given in an appendix, relies on 
an approximation lemma which concerns the behaviour of the Coulomb energy under mollification of measures, 
and on a generalization of an identity of Mattila [Ma95] (see the appendix).

The lemma readily yields
\begin{lemma} \label{convex} $I$ is strictly convex on ${\cal C}(\R^3\backslash\Omega)=
\{\mu\in {\cal M}(\R^3\backslash \Omega) \, | \, 
\mu\ge 0, J(\mu)<\infty\}$.
\end{lemma}
{\bf Proof} Because the first term of $I$ is linear and the second term is quadratic, we have
$$
   \frac{I(\mu_1)+I(\mu_2)}{2} - I\Bigl(\frac{\mu_1 + \mu_2}{2}\Bigr) = \frac14 J(\mu_1-\mu_2)
$$
for any $\mu_1,\mu_2$ in the above set. The assertion now follows from Lemma \ref{positive}.

Finally we will need the following, much simpler, result, which does not rely on Lemma \ref{positive}:
\begin{lemma} \label{lowersemi}
$I$ is (sequentially) weak* lower semicontinuous on ${\cal C}(\R^3\backslash\Omega)$, i.e. if $\mu$, $\mu_j\in{\cal
C}(\R^3\backslash\Omega)$ with 
$\mu_j\weakstar\mu$, then $I(\mu)\le\liminf_{j\to\infty} I(\mu_j)$. 
\end{lemma}
{\bf Proof} This follows, e.g., from Theorem \ref{theo3} and the general fact that $\Gamma$-limits are
lower semicontinuous (see [Br02, Proposition 1.28]). To keep this section self-contained,
we include a direct proof, via a simple truncation argument which
replaces the discontinuous integrand $1/|x-y|$ in $I$ by a continuous function. Let
\be\label{falpha}
    f^\alpha(x,y):=\left\{\begin{array}{lll} 
    \frac{1}{|x-y|}\quad\mbox
 { if } |x-y|\geq \alpha \\
 \quad\frac{1}{\alpha}\quad\mbox{ if } |x-y|\leq\alpha,\qquad
 \end{array}\right. 
\ee
and let $I_{\alpha}$ be the functional obtained by replacing the integrand $1/|x-y|$ in the second term of $I$
by $f^\alpha(x,y)$.
Then $\liminf_{j\to\infty}I(\mu_j)\ge\liminf_{j\to\infty}I_\alpha(\mu_j)
\ge I_\alpha(\mu)$, due to the trivial inequality $1/|x-y|\ge f_\alpha(x,y)$ and
the convergences $\mu_j\weakstar\mu$ and $\mu_j\otimes\mu_j\weakstar\mu\otimes\mu$. 
To finish the proof it suffices to show that $\lim_{\alpha\to 0}I_{\alpha}(\mu) = I(\mu)$. 
If $(\mu\!\otimes\!\mu)(\Diag)>0$ then this is true because both sides are equal to $+\infty$; if $(\mu\!\otimes\!\mu)(\Diag)=0$
then this follows by monotone convergence, because $f_\alpha(x,y)$ is monotonically increasing in $\alpha$ and
tends to $1/|x-y|$ for all $(x,y)\neq \Diag$, and hence for $(\mu\!\otimes\!\mu)$-a.e. $(x,y)$. 
\\[2mm]
{\bf Proof of Proposition \ref{limitmodel}} Existence of a minimizer is immediate from Lemma \ref{lowersemi}: any
minimizing sequence $\mu^{(j)}$ is bounded in ${\cal M}(\R^3\backslash\Omega)$, since $\int d\mu^{(j)})\le\lambda$;
thus there exists a weak* convergent subsequence, by the Banach-Alaoglu theorem; its limit must be a minimizer, by
Lemma \ref{lowersemi}. Uniqueness follows from Lemma \ref{convex}. 

Next we show e). The idea is to first replace the delta functions generating the nuclear potential by spherical
charge distributions so as to make their self-interaction finite instead of infinite; one can then rewrite the energy
by ``completing the square''. So let
$$
  \rho := \sum_{\alpha=1}^M z_\alpha \frac{ H^2\bigl|_{S_\alpha}\bigr. }{4\pi d^2},
$$
with $S_\alpha$ as in Theorem \ref{thm1}. By Newton's theorem that the electrostatic potential exerted by
a radial charge distribution onto a point outside it is the same as that exerted by the same amount of charge placed
at the centre of the sphere,
\be \label{newton}
    \sum_{\alpha=1}^M \frac{-z_\alpha}{|x-R_\alpha|} = - \int_{\R^3} \frac{1}{|x-y|}d\rho(y) 
    \;\;\;\mbox{ for all }x\in\R^3\backslash\Omega.
\ee
It follows that for any $\mu\in{\cal M}_+(\R^3\backslash\Omega)$, 
\begin{eqnarray} 
  \hspace*{-5mm} I(\mu) & \! = \! & - \int\int_{\R^6}\frac{1}{|x-y|}d\rho(y)d\mu(x) + \frac12 \int\int_{\R^6}\frac{1}{|x-y|}d\mu(x)d\mu(y) \nonumber \\
  \hspace*{-5mm}        & \! = \! & \frac12 \int\int_{\R^6} \frac{1}{|x-y|} d(\rho-\mu)(x)d(\rho-\mu)(y) - \frac12 \int\int_{\R^6} 
                \frac{1}{|x-y|}d\rho(x)d\rho(y). \label{squaretrick}
\end{eqnarray}
Since $\int d\rho = z$ and, by assumption, $z\le \lambda$, the measure $\mu=\rho$ is contained in $\Mlambda$; hence
by Lemma \ref{positive} it is the unique minimizer of $I$ on $\Mlambda$. 

It remains to evaluate $I(\rho)$. Using (\ref{squaretrick}) , denoting $\rho_\alpha=z_\alpha (4\pi d^2)^{-1}H^2|_{S_\alpha}$,
and again using Newton's theorem,
\begin{eqnarray*}
   I(\rho) &=& -\frac12\int\int\frac{d\rho(x)d\rho(y)}{|x-y|} \\
   &=& -\frac12 \sum_{\alpha}\int\int\frac{d\rho_\alpha(x)d\rho_\alpha(y)}{|x-y|} - \sum_{\alpha<\beta}\int\int
     \frac{d\rho_\alpha(x)d\rho_\beta(y)}{|x-y|} = -\frac12\sum_\alpha \frac{z_\alpha^2}{d} - \sum_{\alpha<\beta}
     \frac{z_\alpha z_\beta}{|R_\alpha-R_\beta|}.
\end{eqnarray*}
This completes the proof of e), and also establishes d) in case $\lambda\ge z$. 

Next we show d) when $\lambda\le z$. We assume without loss of generality $R_1=0$, $d=1$,
and use a symmetrization argument. Given $\mu\in{\cal C}$, let 
$$
   \bar{\mu} = \frac{1}{\mbox{vol}(SO(3))}\int_{SO(3)} \mu_R d{\cal H}(R),
$$
where ${\cal H}$ is the Haar measure on $SO(3)$ and $\mu_R$ is the rotated measure $\mu_R(A) = \mu(R^{-1}A)$
($R\in SO(3), \, A$ measurable). By the strict convexity of $I$ and Jensen's inequality, 
\be \label{symmetrization}
  I(\bar{\mu}) \le I(\mu), \mbox{ with equality if and only if }\mu=\bar{\mu}.
\ee
Hence it suffices to show that among radially symmetric measures in $\Mlambda$, the unique minimizer of $I$ 
is given by the formula in d). 

But for such $\mu$, we can rewrite $I$ in terms of the radial measure 
$$
   \nu(r) := \mu(S_r), \;\;\; S_r=\{x\in \R^3\, | \, |x|=r\}, 
$$
as follows. (Note that $\mu\in\Mlambda$ is equivalent to $\nu\in \Mlambda^{radial}=\{\nu\in {\cal M}([1,\infty))\, | \, 
\nu\ge 0, \int d\nu\le \lambda\}$, and $\mu\in {\cal M}_\lambda$ is equivalent to $\nu\in
{\cal M}_\lambda^{radial}=\{\nu\in {\cal M}([1,\infty))\, | \, 
\nu\ge 0, \int d\nu = \lambda\}$.)
By Newton's theorem,
$$
   \int_{\{y\, | \, |y|\le |x|\} } \frac{d\mu(y)}{|x-y|} = \Bigl( \int_{[1,|x|]} d\nu(|y|) \Bigr) \cdot \frac{1}{|x|},
$$
and by the fact that the electrostatic potential exerted by a radial charge distribution onto a point inside it is constant,
$$
   \int_{\{y\, | \, |y| > |x|\} } \frac{d\mu(y)}{|x-y|} = \int_{(|x|,\infty)}\frac{1}{|y|} d\nu(|y|).
$$
Consequently, letting $B_1:=\{x\in\R^3\, | \, |x|=1\}$, 
$$
   \int_{\R^3\backslash B_1} \frac{d\mu(y)}{|x-y|} = \int_{[1,\infty)}\!\!\! \min\{\frac{1}{|x|},\frac{1}{|y|}\} d\nu(|y|)
$$
and
\begin{equation} \label{itilde}
  I(\mu) = -\int_{[1,\infty)} \frac{z}{r} d\nu(r) + \frac12 \int\int_{[1,\infty)^2} \!\!\! \min\{\frac1r,\frac{1}{r'}\}
       d\nu(r)d\nu(r') =: \Itilde(\nu).
\end{equation}
Hence to complete the proof of d), it suffices to show that 
\begin{equation} \label{numin}
   \nu=\min\{\lambda,z\}\delta_1
\end{equation}
is a minimizer of $\Itilde$ on $\Mlambda^{radial}$. But for $\nu\in {\cal M}_\lambda^{radial}$, since $1=\int d\nu / \lambda$, 
$$
   \Itilde(\nu) = \frac{1}{2}\int\!\int_{[1,\infty)^2}\Bigl( -\frac{z}{\lambda r} - \frac{z}{\lambda r'} +
\min\{\frac{1}{r},\frac{1}{r'}\}\Bigr) d\nu(r)\, d\nu(r').
$$
Now for $\lambda\le z$ the integrand is minimized pointwise at $(r,r')=(1,1)$, so $\nu(r)=\lambda\delta_1$
is a minimizer on ${\cal M}_\lambda^{radial}$. This completes the proof of d). 

Next, we establish a). For $\lambda\ge z$, this follows from the explicit formulae in d), e). Thus it suffices
to show $\int d\mu_\lambda=\lambda$ for $\lambda\le z$. The inequality ``$\le$'' is trivial. To prove ``$\ge$'', 
suppose $\int d\mu_\lambda < \lambda$. Then the measure $\tilde{\mu}=\mu + \epsilon (4\pi R^2)^{-1}H^2|_{S_R(0)}$,
$S_R(0)=\{x\in\R^3\, | \, |x|=R\}$, lies in $\Mlambda$ for sufficiently small $\epsilon>0$, and as $R\to\infty$,
$\epsilon\to 0$, by multipole expansion
$$
   I(\tilde{\mu})=I(\mu) + \epsilon \Bigl( \frac{\int d\mu - z}{R} + o\Bigl(\frac1R\Bigr)\Bigr) + O(\epsilon^2),
$$
which is smaller than $I(\mu)$ for $R$ sufficiently large and $\epsilon$ sufficiently small, since
$\int d\mu - z < \lambda - z \le 0$. This contradicts the minimality of $\mu$, completing the proof of a).

b) is a straightforward consequence of a): if $\lambda\ge z$ and $\mu_\lambda$ minimizes $I$ on $\Mlambda$, 
then by a) $\mu_\lambda\in{\cal M}_{[0,z]}$, and hence minimizes $I$ on ${\cal M}_{[0,z]}$. Uniqueness now implies
$\mu_\lambda=\mu_z$. 

It remains to prove c). That $e(\lambda)$ is monotonically nonincreasing in $\lambda$ is trivial. 
Continuity at $\lambda$ follows by using the measure  $\mu=(\lambda/\lambda')\mu_{\lambda'}$, $\lambda'>\lambda$, as trial function in the
variational principle `Minimize $I$ on $\Mlambda$' and letting $\lambda'\to\lambda$. That $e(\lambda)$ is constant
for $\lambda\in[z,\infty)$ is obvious from b). Finally, strict monotonicity on $[0,z]$ can be seen as follows. Let
$0\le\lambda<\lambda'\le z$. By monotonicity, $e(\lambda)\ge e(\lambda')$. But ``='' is impossible, since otherwise 
uniqueness of minimizers in ${\cal M}_{[0,\lambda']}$ would enforce $\mu_\lambda=\mu_{\lambda'}$, contradicting the fact
that by a), $\int d\mu_\lambda \neq \int d\mu_{\lambda'}$. 

The proof of the proposition is complete. 
\\[2mm]
While from the point of view of Gamma-convergence it is natural to work in Proposition \ref{limitmodel} above
with the relaxed constrained $0\le\int d\mu\le \lambda$ (see Theorem \ref{theo3}), it is also of interest to consider
the sharp constraint $\int d\mu=\lambda$, for this yields a strikingly simple continuum version of the
attainment/nonattainment transition of the particle system described in Proposition \ref{prop1}a) and Corollary 
\ref{cor1}:
\setcounter{theorem}{0}
\begin{corollary} Assume condition (\ref{nonoverl}) holds. 
Let ${\cal M}_\lambda := \{\mu\in {\cal M}(\R^3\backslash \Omega)\, | \, \mu\ge 0, \int d\mu = \lambda\}$. 
If $\lambda\le z$, then the unique minimizer of $I$ on ${\cal M}_\lambda$ is given by $\mu_\lambda$. 
If $\lambda > z$, then the infimum of $I$ on ${\cal M}_\lambda$ is not attained, the value of the infimum equals
$\min_{{\cal M}_z} I$, and any minimizing sequence $\mu^{(j)}$
converges weak* but not strongly to $\mu_z$; in particular $\int d\mu_z = z < \lambda = \lim_{j\to\infty}\int d\mu^{(j)}$. 
\end{corollary}
The proof is straightforward from Proposition \ref{limitmodel} and Lemma \ref{lowersemi}. See [Ba88]
for a unifying mathematical
setting which subsumes many examples of such attainment/nonattainment transitions associated with
a ``loss of the constraint'' phenomenon.

\section{Neutrality and equidistribution}
Here we show how the explicit solution of the limit theory derived above, combined with the abstract
Gamma-convergence result of Theorem \ref{theo3}, allows infer the abstract convergence result of Theorem \ref{thm2}
and the neutrality and equidistribution
results for the particle system stated in Theorem \ref{thm1}.

The main point is that in the topology 
in which the Gamma-convergence occurs, the sequence of (associated measures of) approximate minimizers
of the particle energy $\vnz$ is compact. The rest of the argumentation is standard in Gamma-convergence.
\\[2mm]
{\bf Proof of Theorem \ref{thm2}} 
First, note that by (\ref{Inzid}) and Theorem \ref{theo3}, 
\be \label{limit'}
    \frac{\inf V_{N,Z}}{Z^2} = \inf_{\Mlambda} \inz \to \inf_{\Mlambda} I
\ee
in the limit (\ref{limit}). 

Now let $(x_1^{(N,Z)},...,x_N^{(N,Z)})$ be a sequence of approximate minimizers of $\vnz$. It follows from
definition (\ref{approxmin}) and (\ref{limit'}) that
\be \label{limit''}
   \inz(\xonenz,...,\xnnz) \to \inf_{{\cal M}_{[0,\lambda]}} I.
\ee
In addition the associated sequence of measures $\munz$ defined in (\ref{measures}) is bounded in 
${\cal M}(\R^3\backslash \Omega)$ (because $\munz\ge 0$ and $\int_{\R^3\backslash \Omega}d\munz = N/Z$ is
bounded), and hence weak* compact, by the Banach-Alaoglu theorem. 

By the lower bound property (ii) contained in the Gamma-convergence result of Theorem \ref{theo3}, for every
weak* convergent subsequence $\munjzj$ the limit $\tilde{\mu}$ satisfies $I(\tilde{\mu})\le\liminf I(\munjzj)$.
But by (\ref{limit''}) and Proposition \ref{limitmodel}, $\tilde{\mu}$ must equal the unique minimizer
of $I$ on ${\cal M}_{[0,\lambda]}$. 

Finally, since every subsequence of $\munz$ converges to this minimizer, so must the whole sequence. 

This establishes the theorem.
\\[2mm]
{\bf Proof of Theorem \ref{thm1}} This is a straightforward consequence of Theorem \ref{thm2}. Specializing to $N=Z$
(whence $\lambda = 1$) and using Proposition \ref{limitmodel} and the fact that $\sum_{\alpha=1}^Mz_\alpha=1$ shows that 
$\munz\weakstar \sum_{\alpha=1}^M z_\alpha(4\pi d^2)^{-1}H^2|_{S_\alpha}$ in ${\cal M}(\R^3\backslash\Omega)$. 
But by Proposition \ref{prop1}, the measures $\munz$ are supported on $\cup_\alpha S_\alpha$, so the
above convergence also occurs in ${\cal M}(\cup_\alpha S_\alpha)$. 
We now use the well known fact that if a sequence of Radon measures $\mu^{(j)}$ on any compact $d$-dimensional manifold $X$ converges
weak* to $\mu$, then for all Borel sets $A\subseteq X$ with $\mu(\partial A)=0$, $\mu^{(j)}(A)\to\mu(A)$.
Consequently for $A\subseteq S_\alpha$
$$
   \frac{ \sharp \{ x_i^{(N)} \, | \, x_i^{(N)}\in A\} }{N} = 
   \int_{S_\alpha} \chi_{_A}(x)d\munz(x) \to \int_{S_\alpha} \chi_{_A}(x) d\mu(x) = \frac{ area(A) }{area(S_\alpha)}.
$$
This proves Theorem \ref{thm1}. 

\section{Instability of asymptotically negative ions}
Here we show how Theorem \ref{thm2} together with the saturation of mass phenomenon of Proposition \ref{limitmodel} a) 
leads to a nonattainment result on $\vnz$ which complements Proposition \ref{prop1}a). 
\setcounter{theorem}{0}
\begin{corollary} \label{cor1}
(Instability of asymptotically negative ions) \\ Let $N_*(\underline{Z}) = \sup\{N\in\N \, | \, \inf_{\an} \vnz
\mbox{ attained}\}$. Then as $Z=\sum_{\alpha=1}^MZ_\alpha\to\infty$, $Z_\alpha/Z\to z_\alpha$, 
$$
     \frac{N_*(\underline{Z})}{Z} \longrightarrow 1.
$$
\end{corollary}
\setcounter{theorem}{2}
{\bf Proof}
The simple attainment result in Proposition \ref{prop1} a) implies that $N_*(\underline{Z})\ge Z+1$; hence it is
clear that $\liminf_{Z\to\infty}\frac{N_*(\underline{Z})}{Z} \ge 1$. The nontrivial assertion in the corollary is that
\be \label{limsup}
   \limsup_{Z\to\infty} \frac{N_*(\underline{Z})}{Z} \le 1. 
\ee
But this is a direct consequence of Theorem \ref{thm2} and Proposition \ref{limitmodel}a), as follows.
Denote the value of the $\limsup$ in
(\ref{limsup}) by $\lambda_*$, and consider any subsequence realizing it, i.e. $N_*(\underline{Z}^{(j})/Z^{(j)}\to\lambda_*$
as $j\to\infty$. Abbreviate $N_j:=N_*(\underline{Z}^{(j})$, and let
$(x_1^{(j)},...,x_{N_j}^{(j)})$ 
be a minimizer of $V_{N_j,\underline{Z}^{(j)}}$. 

On the one hand, by Theorem \ref{thm2}, Proposition \ref{limitmodel}a), and the fact that $\sum_\alpha z_\alpha=1$,
the associated measure satisfies 
$\mu^{(N_j,\underline{Z}^{(j)})}\weakstar\mu$ in ${\cal M}(\R^3\backslash\Omega)$, for some measure
$\mu$ with $\int_{\R^3\backslash\Omega}d\mu = \min\{\lambda_*,1\}$. 

On the other hand, by the absorption principle (Proposition \ref{prop1} b)), letting $S:=\cup_\alpha S_\alpha$,
$\int_{S}d\mu^{(N_j,\underline{Z}^{(j)})}\to \lambda_*$ and, due to the fact that the above weak* convergence also occurs in
${\cal M}(S)$, $\int_S d\mu = \lambda_*$. 

Consequently $\min\{\lambda_*,1\}=\lambda_*$, or equivalently $\lambda_*\le 1$, as was to be shown.
\section{The continuum theory as Gamma limit of the many-particle Coulomb system}
Here we show that the continuum theory (\ref{contlimit}) arises in a mathematically rigorous way
(namely as a Gamma-limit) from the many-particle Coulomb system, i.e. we prove Theorem \ref{theo3}.
As emphasized in the Introduction, despite a result of this kind being --- in our view --- very natural,
we know of no case of any many-atom or many-electron system where such a result has been
previously established.

Recall (e.g. from [Br02]) that a sequence $I^{(j)}:X\rightarrow\R\cup\left\{\infty\right\}$ of functionals on a
topological space $X$ is said to $\Gamma$-converge to 
$I \, : \, X\rightarrow\R\cup\left\{\infty\right\}$ if for all $\mu\in X$ we have:
\\[1mm]
(i) (Ansatz-free lower bound) For every sequence $\mu_j\in X$ converging to $\mu$ we have
$I(\mu)\leq\liminf_{j\rightarrow\infty} I^{(j)}(\mu_j)$.
\\[1mm]
(ii) (Attainment of lower bound) There exists a sequence $\mu_j\in X$ converging to $\mu$ such that
$I(\mu) = \lim_{j\to\infty}I^{(j)}(\mu_j)$.
\\[2mm]
In our case, $X=\Mlambda$ (see Theorem \ref{theo3}), the space of nonnegative Radon measures on $\R^3\backslash \Omega$
of mass $\le \lambda$, endowed with the weak* topology, 
and $I^{(j)} = I^{(N^{(j)},\underline{Z}^{(j)})}$, where $N^{(j)}\in\N$, $\underline{Z}^{(j)}=(Z_1^{(j)},...,Z_M^{(j)})$, 
$Z^{(j)}=\sum_{\alpha=1}^M Z_\alpha^{(j)}$, and 
\be \label{seqlimit}
   \frac{N^{(j)}}{Z^{(j)}}\to\lambda \in(0,\infty), \;\;\; \frac{Z_\alpha^{(j)}}{Z^{(j)}}\to z_\alpha
\;\;\;(j\to\infty).
\ee
In order to establish Theorem \ref{theo3} we need to verify (i) and (ii). 
\subsection{Proof of the lower bound (i)}
Suppose that $\mu_j\weakstar\mu$. We may assume without loss of generality that $I^{(j)}(\mu_j)<\infty$ for all $j$
(because if $J := \{j\in\N\, | \, I^{(j)}(\mu_j)<\infty\}$ is finite the assertion is trivial and if it is infinite
then $\liminf_{j\to\infty} I^{(j)}(\mu_j) = \liminf_{j\in J, \, j\to\infty}I^{(j)}(\mu_j)$). By passing to a subsequence
we may in addition assume that $I^{(j)}(\mu_j)\to \liminf_{j\to\infty}I^{(j)}(\mu_j)$. 

We use a truncation argument as in the proof of Lemma \ref{lowersemi}. We let $f_\alpha$, $I_\alpha$ be as defined there, and compute
\begin{eqnarray*}
  I^{(k)}(\mu_k) & = & - \int_{\R^3\backslash \Omega}
\sum_{\alpha=1}^M\frac{Z_\alpha^{(k)}}{Z^{(k)}|x-R_\alpha|}d\mu_k(x) +\frac12 
  \int\int_{(\R^3\backslash \Omega)^2\backslash \Diag} \frac{1}{|x-y|}d\mu_k(x)d\mu_k(y) \\
  & \ge & 
  - \int_{\R^3\backslash \Omega}  \sum_{\alpha=1}^M\frac{Z_\alpha^{(k)}}{Z^{(k)}|x-R_\alpha|}d\mu_k(x) +\frac12 
  \int\int_{(\R^3\backslash \Omega)^2\backslash \Diag} f^\alpha(x,y)d\mu_k(x)d\mu_k(y) \\
  & = & I_{\alpha}(\mu_k)  + \sum_{\alpha=1}^M \Bigl(z_\alpha - \frac{Z_\alpha^{(k)}}{Z^{(k)}}\Bigr)\int_{\R^3\backslash\Omega}
   \frac{1}{|x-R_\alpha|}d\mu_k(x) - \frac{N^{(k)}}{2(Z^{(k)})^2\alpha}.
\end{eqnarray*}
Using the fact that the last term on the right hand side and the factors in the middle term tend to zero by (\ref{limit}) and
that if $\mu_k\weakstar\mu$ in ${\cal M}(\R^3\backslash \Omega)$ then 
$\mu_k\otimes\mu_k\weakstar\mu\otimes\mu$ in ${\cal M}((\R^3\backslash \Omega)^2)$, letting $k$ tend to infinity gives 
\be
 \lim_{k}I^{(k)}(\mu_k)\ge
 \lim_{k}I_{\alpha}(\mu_k)=I_{\alpha}(\mu).
\ee
But as shown at the end of the proof of Lemma \ref{lowersemi}, 
$\lim_{\alpha\to 0}I_{\alpha}(\mu) = I(\mu)$. 
This establishes (i).
\subsection{Proof of the upper bound (ii)}  
Fix a sequence $N^{(\nu)}\to\infty$, $\underline{Z}^{(\nu)}$ satisfying (\ref{seqlimit}).
Given $\mu\in {\cal M}_+(\R^3\backslash\Omega)$, we need to construct a sequence $\mu_\nu\in
{\cal M}_+(\R^3\backslash \Omega)$ (or, in Gamma-convergence terminology, a ``recovery sequence'')
such that $\mu_\nu\weakstar\mu$ and
$\limsup_{\nu\to\infty} \inznu(\mu_\nu)\le I(\mu)$. 

This is achieved by a careful multiscale construction, by introducing a mesoscale $h$ with
particle spacing $<<$ mesoscale $<<$ diameter of support of $\mu$ (see Step 2) and approximating $\mu$ in each mesoscale
region by a suitable number of Dirac masses placed on some suitable lattice (see Step 4).
The number of Dirac masses is governed by the amount of mass to be accommodated in the region
(see Step 3).

%

A first difficulty is that unlike in usual arguments establishing density of discrete measures, the amplitude of each
Dirac mass is fixed exactly to be $1/Z$, which leads to a mass error of order $1/Z$ in a typical mesoscale region. 

A second, and more fundamental, difficulty is that
one expects the energy to be highly sensitive to the precise placement of the
particles; but the precise structure of approximate or exact minimizers
of the many-body Coulomb interaction is unknown mathematically.
We know of no
attempt to prove minimizers are cystalline, or approximately crystalline, let alone to establish the optimal
lattice structure -- neither for Coulomb interactions nor for any other realistic interaction law in three
dimensions.

Very remarkably, the long range nature of the Coulomb force, usually considered a complicating rather than
a simplifying feature, works in our favour. It implies that the energy is dominated by long range contributions,
and so at short range a rough knowledge of bondlengths (to within a factor) turns out to be sufficient, 
as long as knowledge of the long range distances, governed by the ``packing density'', is precise. The key
point in the proof is the implementation of these ideas in Step 5.
\\[2mm]
{\bf Step 1} {\it Reduction to compactly supported measures of finite energy, mass $\lambda$, and bounded Lebesgue
density} \\[1mm]
If $I(\mu)=\infty$ existence of a recovery sequence is trivial: for instance the sequence $\mu_\nu\equiv\mu$
will do. So we may assume $I(\mu)<\infty$. 

By a standard approximation argument, we may also assume that $\mu$ has compact support in
$\R^3\backslash\overline{\Omega}$.

A little less trivially, we claim that it is enough to establish existence of a recovery sequence for
measures with $\int d\mu=\lambda$. This is because if $\int d\mu>\lambda$, then $I(\mu)=\infty$
and we are back in the case dealt with above, whereas if $\int d\mu <\lambda$ we can always ``place
unwanted mass at infinity''. More precisely, if $\int d\mu =\tilde{\lambda} <\lambda$, we choose
$\tilde{N}_\nu\in\N$, $\tilde{N}_\nu<N_\nu$, such that $\tilde{N}_\nu/Z_\nu\to\tilde{\lambda}$, apply 
existence of a recovery sequence $\tilde{\mu}_\nu$ for 
$\mu$ with respect to the functionals $I^{(\tilde{N}_\nu,Z_\nu)}$, and set 
$\mu_\nu = \tilde{\mu}_\nu + \sum_{k=1}^{N_\nu - \tilde{N}_\nu}(1/Z_\nu)\delta_{x_\nu^k}$ with
$\min_k|x_\nu^k|\to\infty$, $\min_{k\neq\ell}|x_\nu^k-x_\nu^\ell|\to\infty$. It follows that
$\mu_\nu-\tilde{\mu}_\nu \weakstar 0$ and $\inznu(\mu_\nu)-I^{(\tilde{N}_\nu,Z_\nu)}(\tilde{\mu}_\nu)\to 0$. 

Finally we claim that we may assume that $\mu$ has bounded Lebesgue density, i.e. $d\mu(x)=m(x)dx$ for some
$m\in L^\infty(\R^3\backslash\Omega)$. (We thank Stefan M\"uller for this idea,
which facilitates the simple energy error estimate (\ref{smestimate})
via the uniform bound (\ref{abound}) on the local lattice spacing below. It 
replaces our original more complicated energy estimate via an integral bound on the lattice spacing.)
This is because given any measure 
$\mu\in {\cal M}_+(\R^3\backslash\Omega$ of mass $\lambda$ and compact support,
one can construct a sequence of measures $\mu^{(\epsilon)}\in
{\cal M}_+(\R^3\backslash\Omega)\cap L^\infty(\R^3)$ of mass $\lambda$ and compact support with $\mu^{(\epsilon)}\weakstar\mu$
and
$\liminf_{\epsilon\to 0}I(\mu^{(\epsilon)})\le I(\mu)$. Indeed, we claim that the mollified measure
$$
   \mu^{(\epsilon)} := \phi_\epsilon * \mu = \int_{B_\epsilon(0)}\phi_\epsilon(z)\mu(\cdot - z)dz
$$
with $\phi_\epsilon(z)=\epsilon^{-3}\phi(\epsilon^{-1}z)$, $\phi\in C_0^\infty(B_1(0))$, $\phi\ge 0$,
$\int_{\R^3}\phi = 1$, has the required properties. For the elementary proof of weak* convergence to $\mu$
see e.g. [M] Theorem 1.26. To verify that $\liminf_{\epsilon\to 0}I(\mu^{(\epsilon)})\le I(\mu)$, we write
$I(\mu) = \int_{\R^3\backslash\Omega}v(x)d\mu(x) + J(\mu)$, with $J$ as in (\ref{Jfunctional}). The weak* 
convergence of $\mu^{(\epsilon)}$ implies $\int_{\R^3\backslash\Omega}v(x)d\mu^{(\epsilon)}(x)\to
\int_{\R^3\backslash\Omega}v(x)d\mu(x)$, and the translation invariance and convexity of $J$ and Jensen's
inequality imply 
$$
   J(\mu) = \int_{B_\epsilon(0)}\phi_\epsilon(z)J(\mu(\cdot - z))\, dz \ge J\Bigl(
            \int_{B_\epsilon(0)}\phi_\epsilon(z)\mu(\cdot - z)\, dz\Bigr) = J(\mu^{(\epsilon)})
$$
for all $\epsilon$, establishing the assertion. We remark that by weak* lower semi-continuity of $I$, in fact
one has $I(\mu^{(\epsilon)})\to I(\mu)$, but this is not needed here.
\\[2mm]
{\bf Step 2} {\it Discretization of $\R^3\backslash\Omega$ into mesh of size $h$} \\[1mm]
From now on we fix a measure $\mu\in{\cal M}_+(\R^3\backslash\Omega)$ which has finite energy and
is compactly supported in $\R^3\backslash\Omega$. Hence for $h_0$ sufficiently small, supp $\mu$ is contained 
in some finite union of disjoint cubes of sidelength $h_0$, $\Omega'=\cup_{\alpha=1}^{N_0}Q_{h_0}(R_\alpha')$, 
with $\Omega'\subset\R^3\backslash\Omega$,
where $Q_h(R)$ denotes a cube centred at $R$ of sidelength $h$,  
$\{x\in\R^3\, | \, -h/2 \le (x-R)\cdot e_j <h/2\}$, and where $e_1$, $e_2$, $e_3$ are the standard
basis vectors of $\R^3$. 

Now given any $n\in\N$ (to be chosen later, depending on $N$ and $Z$)
we obtain a mesh of size $h:=h_0/n$ by dividing each cube of sidelength $h_0$ into $n^3$ smaller cubes of sidelength
$h_0/n$. This way we obtain a disjoint family of cubes $\{\Qi\}_{i=1}^{N_0 (h/h_0)^{-3}}$ of sidelength $h$ whose 
union contains $\supp\mu$.
\\[2mm]
{\bf Step 3} {\it Choice of number of Dirac masses in each region and mass error analysis}\\[1mm]
For given $\mu$, $N$, $Z$, we need to approximate $\mu|_{\Qi}$ by a measure of form
\be \label{ansatz}
    \munz\Bigl|_{\Qi}\Bigr. = \frac1Z \sum_{k=1}^{L_i} \delta_{\xik}, \;\;\; L_i\in\N\cup\{0\},\;
    \sum_i L_i=N,
\ee
because otherwise $\inz(\munz)$ is infinite. In particular, the allowed mass has to be an integer multiple of
$1/Z$, enforcing a mass error. Here we deal with the latter, postponing the choice of positions $\xik$
to the next step.

It will be convenient to introduce a trivial amplitude factor
\be\label{amplitudefactor}
   \phi(N,Z) := \frac{N}{Z} \cdot \lambda^{-1}
\ee
(note $\phi(N,Z)\to 1$ in the limit (\ref{limit})) and approximate the measure $\phi\cdot\mu$, because 
$\int d(\phi\mu)=N/Z=\int d\munz$. 

Choose $\ell_i\in\N\cup\{0\}$ such that $\ell_i/Z$ is a good approximation to the mass of $\phi\mu$ 
in $\Qi$, i.e.
\be \label{massbound}
    \frac{\ell_i}{Z} \le \phi\mu(\Qi) \le \frac{\ell_i + 1}{Z}.
\ee
This implies $\sum_i\ell_i/Z \le N/Z \le \sum_{i\, | \, \mu(\Qi)>0}(\ell_i+1)/Z$. Hence $(\sum_i\ell_i) + r=N$
for some integer $r$ less or equal to the number of indices $i$ with $\mu(\Qi)>0$. It follows that if for
$r$ such indices we set $L_i:=\ell_i+1$, and let $L_i:=\ell_i$ otherwise, we have $\sum_iL_i=N$, 
$L_i=0$ when $\phi\mu(\Qi)=0$, and, by (\ref{massbound}), 
\be \label{Liest}
   \frac{L_i-1}{Z} \le \phi\mu(\Qi) \le \frac{L_i+1}{Z}.
\ee
Hence if $\munz|_{\Qi}$ is given by (\ref{ansatz}), regardless of the choice of positions $\xik$
we have
\be \label{constbound}
    \Bigl| \int_{\Qi} A\, d(\phi\mu) - \int_{\Qi} A\, d\munz \Bigr| \le \frac{|A|}{Z} \;\;\mbox{ for all }
    A\in\R.
\ee
Let $c_i$ denote the centre of the cube $\Qi$. Then by (\ref{constbound}),
if $g\in C(\R^3\backslash\Omega)$ and $\delta$ is its modulus of continuity on the lengthscale of the mesh, 
$$
   \delta := \sup_{|x-y|\le \max_i\diam\Qi} \Bigl| g(x)-g(y) \Bigr| ,
$$
then 
\begin{eqnarray} \label{contest}
\hspace*{-8mm}  & &  \Bigl| \int_{\R^3\backslash\Omega} g (\phi d\mu) - \int_{\R^3\backslash\Omega} g\, d\munz \Bigr| =  \nonumber \\ 
\hspace*{-8mm}  & & 
  \Bigl| \sum_i\int_{\Qi}\!\Bigl( (g(x)\!-\! g({c}_i)\Bigr)d(\phi\mu)(x) \,
   + \, g({c}_i)\Bigl(d(\phi\mu)(x) \!-\! d\munz(x)\Bigr) \, + \, \Bigl(g({c}_i)\!-\! g(x)\Bigr)d\munz(x) \Bigr)\Bigr| 
     \nonumber \\
\hspace*{-8mm}  & & \le \delta \int d(\phi\mu) \;\; + \;\; \sum_i\frac{|\sup|g|}{Z} \;\; + \;\; \delta\int d\munz 
  \;\; = \;\; 2\delta \frac{N}{Z} + \sup|g|\frac{N_0h_0^3}{h^3 Z}
\end{eqnarray}
where the factor $N_0h_0^3/h^3=N_0n^3$ in the last term is the number of cells $\Qi$.

If $N\to\infty$, $Z\to\infty$, $h=h(N,Z)$, it follows that the right hand side tends to zero for all
$g\in C_0(\R^3\backslash\Omega)$ {\it provided} the meshsize $h$ satisfies 
\be \label{i}
     h\to 0 \;\;\; \mbox{(``small mesh'')}
\ee
(whence $\delta\to 0$), 
\be \label{ii}
    \frac{Z^{-1/3}}{h}\to 0 \;\;\; \mbox{(``particle spacing smaller than mesh'')}
\ee
(whence the second term in (\ref{contest}) tends to zero). Hence if (\ref{i}) and (\ref{ii}) hold, then
$\munz-\phi(N,Z)\mu\weakstar 0$ and hence, thanks to $\phi(N,Z)\to 1$, $\munz\weakstar\mu$.

To understand the meaning of (\ref{ii}), it is instructive to consider the case when $\mu$ is the uniform measure on 
some region of finite diameter, and when $\munz$ is positioned on a periodic lattice in this region. 
Because $\munz$ has $N$ Dirac masses, the lattice spacing must be $\sim N^{-1/3} \sim Z^{-1/3}$ and hence
$Z^{-1/3}/h\sim$ (particle spacing)/meshsize. 
\\[2mm]
{\bf Step 4} {\it Choice of positions of Dirac masses} \\[1mm]
The Dirac masses in each $\Qi$ will be positioned by placing a lattice of spacing $\sim (\max_iL_i)^{-1/3}$ in $\Qi$. First, we estimate $\max_iL_i$, using the fact that $\mu$ has
bounded Lebesgue density. By (\ref{Liest}) and the fact that
$\mu(\Qi)\le h^3||\mu||_{L^\infty}$, 
$$
   L_i \le Z\phi(N,Z) ||\mu||_{L^\infty} h^3 + 1.
$$
Hence since $h^3Z\to\infty$ and $\phi(N,Z)\to 1$, we may assume
\be\label{Lplusdef}
   L_i\le Ch^3 Z =:L_+, \;\;\; C=2||\mu||_{L^\infty}.
\ee
Now we choose a lattice $a\Z^3$ of sufficiently small spacing $a$ so that each $\Qi$ contains at least $L_i$ lattice points.
Because $\Qi$ has sidelength $h$, it suffices to take 
\be \label{adef}
   a = \frac{h}{\lceil(L_+)^{1/3}\rceil},
\ee
where 
$\lceil (L_+)^{1/3}\rceil$ denote the smallest integer $\ge (L_+)^{1/3}$. For future use we note that since $L_+\ge 1$,
$\lceil (L_+)^{1/3} \rceil \le 2 (L_+)^{1/3}$ and so
\be \label{abound}
    \frac{1}{2C^{1/3}Z^{1/3}} = \frac{h}{2L_+^{1/3}} \le a \le \frac{h}{L_+^{1/3}} = \frac{1}{C^{1/3}Z^{1/3}}.
\ee
Now for each $\Qi$, we choose $\{ x_i^{(1)},...,x_i^{(L_i)}\}$ as a subset of $a\Z^3\cap \Qi$,  
and let $\munz$ be as defined in (\ref{ansatz}).  
\\[2mm]
{\bf Step 5} {\it Analysis of energy error}
\\[1mm]
With the above choice of the $\xik$, with $h=h_\nu$ chosen to satisfy (\ref{i}) and (\ref{ii}), 
and with $a=a_\nu$ as defined by (\ref{adef}), (\ref{Lplusdef}), we claim that 
\be \label{mainclaim}
   \lim_{\nu\to\infty} \inznu(\munznu) = I(\mu).
\ee
By the weak* convergence of $\munznu$ to $\mu$, we can immediately pass to the limit in the electron-nuclei
interaction:
$$
    \int_{\R^3\backslash\Omega} v(x)\, d\munznu(x) \;\; \longrightarrow \;\; 
    \int_{\R^3\backslash\Omega} v(x)\, d\mu(x).
$$
As for the electron-electron interaction, we decompose
\begin{eqnarray*}
 & &  \frac12\int\int_{(\R^3\backslash\Omega)^2\backslash\Diag} \frac{1}{|x-y|}d\munznu(x)d\munznu(y) 
 - \frac12\int\int_{(\R^3\backslash\Omega)^2\backslash\Diag} \frac{1}{|x-y|}d\mu(x)d\mu(y) \\
 & & \hspace*{-1cm} =  \frac12 \int\int_{(\R^3\backslash\Omega)^2} f_\alpha \Bigl( d\munznu\otimes d\munznu 
     - d\mu\otimes d\mu \Bigr) 
 \;\;\; - \;\;\; \frac{1}{2\alpha}\frac{N_\nu}{(Z_\nu)^2} \\
 & & \hspace*{-1cm} + \frac12 \int\int_{(\R^3\backslash\Omega)^2\backslash\Diag} \Bigl( \frac{1}{|x-y|}-f_\alpha\Bigr) 
                                                                  d\munznu\otimes d\munznu
  - \frac12 \int\int_{(\R^3\backslash\Omega)^2}  \Bigl( \frac{1}{|x-y|}-f_\alpha\Bigr) d\mu\otimes d\mu.
\end{eqnarray*}
Because the measures $\munznu$ are supported in $\Omega'$, and $f_\alpha$ is continuous
on the closure of $\Omega'\times\Omega'$,
it follows from the weak* convergence of $\munznu$ that the first term tends to zero as $\nu\to\infty$. 
The second and fourth term are $\le 0$, and the integrand in the third term is bounded from above by
$\frac{1}{|x-y|}\chi_{\{|x-y|<\alpha\} }$. It follows that
\be \label{boundbymuj}
   \limsup_{\nu\to\infty} \Bigl(I^{(N_\nu,Z_\nu)}(\munznu) - I(\mu)\Bigr) \le 
   \limsup_{\nu\to\infty} \frac12 \int\!\int_{{(\R^3\backslash\Omega)^2\backslash\Diag}\atop{
    \cap \{|x-y|<\alpha\} } }\!
       \frac{1}{|x-y|}d\munznu \otimes d\munznu 
\ee
for all $\alpha>0$.

It remains to show that the right hand side of (\ref{boundbymuj}) tends to zero as $\alpha\to 0$. By the trivial
estimate for any nonnegative Radon measure $m$
$$
   \int\!\int_{{(\R^3\backslash\Omega)^2\backslash\Diag}\atop{
    \cap \{|x-y|<\alpha\} } }\!
       \frac{1}{|x-y|}dm \otimes dm \le \Bigl(\int_{\R^3}dm\Bigr) \; \cdot \; \sup_{x\in \mbox{supp}\,m}
       \int_{y\in B_\alpha(x)\backslash\{x\} }\frac{1}{|x-y|}dm(y)
$$
and the fact that $\int d\munznu\to\lambda$, it suffices to show that 
\be \label{boundbysingleint}
   \limsup_{\nu\to\infty} \sup_{{x}\in\mbox{supp}\, \munznu} \int_{{y}\in B_\alpha({x})\backslash\{{x}\} }
   \frac{1}{|{x}-{y}|}d\munznu({y})\to 0 \;\;\;(\alpha\to 0).
\ee
We estimate the right hand side by neglecting the fact that only a subset of $a_\nu\Z^3$ carries Dirac
masses and applying the lower bound (\ref{abound})
on $a_\nu$:
\begin{eqnarray} 
   \sup_{{x} \in\mbox{supp}\, \munznu} \int_{{y}\in B_\alpha({x})\backslash\{{x}\} }
   \frac{1}{|{x}-{y}|}d\munznu({y})  
   & \le & Z_\nu^{-1} \sup_{x\in {a_\nu}\Z^3} \!\!\sum_{ 
                       {y\in {a_\nu}\Z^3\backslash\{x\} }\atop{\cap \{|y-x|<\alpha\} }} \!\!\frac{1}{|x-y|} 
   \nonumber \\
   & \le & 8 C \cdot a_\nu^3 \!\!\sum_{y\in a_\nu\Z^3\backslash\{0\}, \, |y|<\alpha} \!\!\frac{1}{|y|}.
   \label{smestimate}
\end{eqnarray}
Now the latter is a Riemann sum and since its ``meshsize'' $a_\nu\le (CZ_\nu)^{-1/3}\to 0$ by (\ref{abound}), 
it converges in the limit $\nu\to\infty$ to $8C\int_{|y|\le \alpha}|y|^{-1}dy$. But explicit evaluation of
this last integral shows that it
tends to zero as $\alpha\to 0$. This establishes (\ref{boundbysingleint}), and completes the proof that the constructed
multiscale lattice measures $\munznu$ constitute a recovery sequence. 

The proof of the Gamma-convergence result Theorem \ref{theo3} is complete. 
\\[5mm]
{\large\bf Appendix: The Coulomb norm on Radon measures without a sign} 
\\[2mm]
We show here that the Coulomb energy on Radon measures without a sign is strictly positive definite,
and hence gives rise to a norm. This was used to determine the minimizer of the continuum energy
in Proposition 2.1. The result is standard ``folklore'' in potential theory (see e.g.[La72]), 
trivial for measures which are either smooth and rapidly decaying or
nonnegative, and surely well known to experts, but we were unable to find a reference.

Our proof is based on two lemmas. The first is
an approximation lemma for Radon measures concerning the behaviour of the Coulomb norm
under mollification; analogous statements are very well known
in Sobolev spaces. The second is a representation formula via the Fourier transform
which is a modest generalization of a corresponding identity in [Ma95, Ch.12] for nonnegative Radon measures 
with compact support;
our argument is different as the restrictions on sign and support are important for the argument in [Ma95].

Let ${\cal M}_+(\R^3)$ denote the set of nonnegative Radon measures of finite mass on 
$\R^3$, and
let ${\cal C}(\R^3)$ denote the subset of such measures $\mu$ whose Coulomb energy $J(\mu)$
defined in (\ref{Jfunctional})
is finite. Define an extension of $J$ to measures with both negative and positive
part, as follows: if $\mu=\mu_1-\mu_2$ with $\mu_1$, $\mu_2\in{\cal C}(R^3)$, set
\begin{eqnarray*}
  J(\mu_1-\mu_2) & := & \frac12 \int\int_{\R^6}\frac{1}{|x-y|}d\mu_1(x)d\mu_1(y)
                    -  \int\int_{\R^6}\frac{1}{|x-y|}d\mu_1(x)d\mu_2(y) \\
                 &   + & \frac12\int\int_{\R^6}\frac{1}{|x-y|}d\mu_2(x)d\mu_2(y).
\end{eqnarray*}
Since the first and last term are finite by assumption, and the integrand in the middle term is nonnegative,
this is well defined as an element of $\R\cup\{-\infty\}$. Denote by $C_0(\R^3)$
the function space
$\{u \, : \, \R^3\to\R \, | \, u\mbox{ continuous}, \, u(x)\to 0 \mbox{ as }|x|\to\infty\}$.
\\[2mm]
{\bf Lemma A1} {\it Let $\phi\in C_0(\R^3)$ be nonnegative, 
radially symmetric, and satisfy $\int_{\R^3}\phi = 1$.
Let $\mu_1$, $\mu_2\in {\cal C}(R^3)$, let $\mu=\mu_1-\mu_2$, and for $\varepsilon>0$ 
let $\mu^\varepsilon$ denote the mollified measure
\be \label{mollif}
    \mu^\varepsilon(x) = (\phi_\varepsilon * \mu)(x)=\int_{\R^3}\phi_\epsilon(x-x')d\mu(x'), 
\ee
where $\phi_\varepsilon(z)=\varepsilon^{-3}\phi(\varepsilon^{-1}z)$.
Then $J(\mu)=\lim_{\varepsilon\to 0} J(\mu^\varepsilon)$.}
\\[2mm]
{\bf Proof} It suffices to consider the middle term in the definition of $J$, i.e. to show
\be \label{claim}
   \int\int_{\R^6} \frac{1}{|x-y|} d\mu_1(x)d\mu_2(y) = \lim_{\varepsilon\to 0}
   \int\int_{\R^6} \frac{1}{|x-y|} d\mu_1^\varepsilon(x)d\mu_2^\varepsilon(y)
\ee
(for the other terms, set $\mu_1=\mu_2$). By Fubini's theorem, 
\be \label{fubini}
   \int\int_{\R^6} \frac{1}{|x-y|} d\mu_1^\varepsilon(x) d\mu_2^\varepsilon(y) 
   = \int\int_{\R^6} \Bigl(
       \int\int_{\R^6}
       \frac{\phi_\varepsilon(x-x')\phi_\varepsilon(y-y')}{|x-y|}dx\, dy \Bigr)
       d\mu_1(x') d\mu_2(y'). 
\ee
Clearly the term in brackets convergers pointwise to $1/|x'-y'|$ as $\varepsilon\to 0$. Moreover setting
$z=x-x'$ and using the radial symmetry of $\phi_\varepsilon$ together with Newton's theorem
that the electrostatic potential exerted by
a spherical charge distribution onto a point inside it is constant, while that exerted
onto a point outside it is the same as that exerted by the same amount of charge placed
at the centre, 
\begin{eqnarray} 
   \int_{\R^3} \frac{\phi_\varepsilon(x-x')}{|x-y|}dx & = & \int_{\R^3} \frac{\phi_\varepsilon(z)}{|z-(y-x')|}dz 
   = \int_{\R^3} \min\{\frac{1}{|z|}, \, \frac{1}{|y-x'|}\} \phi_\varepsilon(z)dz \nonumber \\
   & \le &
   \frac{1}{|y-x'|}\int_{\R^3}\phi_\varepsilon(z)dz = \frac{1}{|y-x'|}. \label{newtonbound}
\end{eqnarray}
Multiplying by $\phi_\varepsilon(y-y')$, integrating over $y$ and applying (\ref{newtonbound}) again yields
$$
   \int\int_{\R^6} \frac{\phi_\varepsilon(x-x')\phi_\varepsilon(y-y')}{|x-y|}dx\, dy \le
   \int_{\R^3}\frac{1}{|y-x'|} \phi_\varepsilon(y-y')dy \le \frac{1}{|x'-y'|}.
$$
Hence by dominated convergence, the right hand side of (\ref{fubini}) 
converges to the left hand side of (\ref{claim}). This establishes the lemma. 

Next, we derive an expression for the Coulomb energy in terms of the Fourier transform, defined for any
nonnegative Radon measure of finite mass, $\mu\in{\cal M}_+(\R^3)$, by 
$$
  \widehat{\mu}(k) := \int_{\R^3}e^{-ik\cdot x}d\mu(x).
$$
Note that for any such measure, its Fourier transform is a bounded continuous function.
\\[2mm]
{\bf Lemma A2} 
{Let $\mu_1$, $\mu_2\in {\cal C}(R^3)$. Then} 
\be \label{fourierid}
      J(\mu_1-\mu_2) = \frac{1}{2} \frac{1}{(2\pi)^3}\int_{\R^3} \frac{4\pi}{|k|^2}|\widehat{\mu_1} - \widehat{\mu_2}|^2 dk.
\ee
{\bf Proof} We approximate $\mu_1$, $\mu_2$ by smoother measures,
as follows. For $\varepsilon>0$, let $\phieps(x):=(2\pi\varepsilon)^{-3/2}e^{-x^2/(2\varepsilon)}$. 
Given any $\mu\in{\cal C}(R^3)$, denote the associated mollified measure (\ref{mollif}) by $\mueps$. 
Then $\mueps\in L^1(\R^3)\cap L^\infty(\R^3)$, because $\int \phieps*\mu = (\int\phieps)\cdot (\int d\mu) = \int d\mu <\infty$
and $\sup (\phieps*\mu) \le (\sup\phieps) \cdot \int d\mu < \infty$. In particular $\mueps\in L^2(\R^3)$, and
so its Fourier transform $\widehat{\mueps}$ is well-defined as an element of $L^2(\R^3)$. By standard Fourier
calculus and the fact that $f(x)=1/|x|$ has Fourier transform $4\pi/|k|^2$, 
\begin{eqnarray}
   & & \widehat\mueps = \widehat\phieps \cdot \hat{\mu}, \nonumber \\
   & & J(\mu_1^{\varepsilon}-\mu_2^{\varepsilon}) = \frac12 \frac{1}{(2\pi)^3} \int_{\R^3} \frac{4\pi}{|k|^2} 
       |\widehat{\phieps}(k)|^2
       |\widehat{\mu_1}(k)-\widehat{\mu_2}(k)|^2 dk. \label{smooth}
\end{eqnarray}
As $\varepsilon\to 0$, because $\widehat{\phieps}(k)=e^{-\varepsilon k^2/2}$ we have that
$|\widehat{\phieps}(k)|^2$ tends monotonically to $1$, and hence by monotone convergence the right hand side of 
(\ref{smooth}) tends to the right hand side of (\ref{fourierid}). On the other hand, by Lemma A1 the left hand side
of (\ref{smooth}) tends to the left hand side of (\ref{fourierid}). This establishes the lemma.

Finally we assert:
\\[2mm]
{\bf Lemma A3} 
{\it Let $\mu_1$, $\mu_2\in {\cal C}(R^3)$. Then $J(\mu_1-\mu_2)\ge 0$, with equality if and only if $\mu_1=\mu_2$.}
\\[2mm]
{\bf Proof} Nonnegativity is clear from (\ref{fourierid}). Moreover the right hand side of (\ref{fourierid}) 
is strictly positive unless $\hat{\mu_1}=\hat{\mu_2}$ Lebesgue-almost everywhere. 
But by continuity of the $\hat{\mu_i}$, this means $\hat{\mu_1}=\hat{\mu_2}$, and hence $\mu_1=\mu_2$. 
The proof is complete.
\\[10mm]
{\bf \large Acknowledgements} The work of S.C. was supported by a Warwick University graduate research fellowship. 
It is a pleasure to thank
James Robinson for helpful discussions, and Stefan M\"uller for a 
simplification of our original construction of a recovery sequence
in Section 5.
\\[10mm]
{\bf \large References}
\\[2mm]
[Ba88] J.M.Ball, Loss of the constraint in convex variational problems, 
Analyse Mathematique et Applications, Gauthier-Villars, 39-53, 1988
\\[2mm]
[BB53] L.S.Bartell, L.O.Brockway, Phys. Rev. {\bf 90}, 833, 1953
\\[2mm]
[Br02] A.Braides, Gamma-convergence for beginners, Oxford lecture series in
mathematics and its applications Vol. 22, 2002
\\[2mm]
[CLL98] I.Catto, C.Le Bris, P.-L. Lions, Mathematical theory of thermodynamic limits, Oxford University Press, 1998
\\[2mm]
[DL67] F.J.Dyson, A.Lenard, Stability of matter. I, J. Math. Phys. (N.Y.) {\bf 8}, 423, 1967
\\[2mm]
[DM88] G.Dal Maso, An introduction to $\Gamma$-convergence, Progress in Nonlinear Differential Equations and Their Applications 8,
Birkh\"auser, Boston, 1988
\\[2mm]
[EO00] \c{S}. Erko\c{c}, H. Oymak, Rules for distribution of pointcharges on a conducting disk, Rapid communication,
Physical Rev E, Vol 62, No.3, R3075-3076, 2000
\\[2mm]
[Fe85] C.Fefferman, The thermodynamic limit for a crystal, Commun. Math. Phys. 98, 289-311, 1985
\\[2mm]
[Fr03] G.Friesecke, The multiconfiguration equations for atoms and molecules: charge quantization and existence of
solutions, Arch. Rat. Mech. Analysis 169, 35-71, 2003
\\[2mm]
[ILS96] I. Iantchenko, E.H.Lieb, H.Siedentop. Proof of a conjecture about atomic and molecular cores related
to Scott's correction, J. Reine Angew. Math. 472, 177-195, 1996
\\[2mm]
[KS98] A.B.Kuijlaars, E.B.Saff, Asymptotics for minimal discrete energy on the sphere, Transactions of the AMS 350 no. 2, 523-538,
1998
\\[2mm]
[La72] N.S.Landkof, Foundations of modern potential theory, Springer Grundlehren der mathematischen Wissenschaften 180, 
Springer-Verlag, 1972
\\[2mm]
[LD68] A.Lenard, F.J.Dyson, Stability of matter. II, J. Math. Phys. (N.Y.) {\bf 9}, 698, 1968
\\[2mm]
[LS77] E.H.Lieb, B.Simon, The Thomas-Fermi theory of atoms, molecules and solids, Adv. Math. 23, 22-116, 1977
\\[2mm]
[LSST] E.H.Lieb, I.M.Sigal, B.Simon, W.Thirring, Approximate neutrality of 
large-$Z$ ions, Commun. Math. Phys. 116, 635-644, 1988
\\[2mm]
[LT75] E.H.Lieb, W.Thirring, Bound for the kinetic energy of fermions which proves the stability of matter,
Phys.Rev.Lett. 35 No. 11, 687-689, 1975
\\[2mm]
[Ma95] P.Mattila, Geometry of sets and measures in euclidean spaces, Cambridge
studies in advanced mathematics Vol. 44, 1995
\\[2mm]
[MDH96] J. R. Morris, D. M. Deaven, K.M. Ho, Genetic algorithm energy minimization for point charges on a sphere, 
Phys Rev B: Condensed matter Vol 53, No.4, 1740-1743
\\[2mm]
[Zh60] G. M. Zhislin, Issledowanie spectra operatora Schredingera dlja sistemy mno\-gih chastits 
[An investigation of the spectrum of a Schr\"odinger operator for a
many-particle system], Trudy Moskovskogo Matematicheskogo Obshchestva 9, 81-120, 1960

\newpage\noindent

\end{document}